\documentclass[11pt]{article}

\usepackage{amsfonts,amsmath,latexsym,color,epsfig}
\newtheorem{theorem}{Theorem}

\input{epsf}

\makeatletter
\def\Ddots{\mathinner{\mkern1mu\raise\p@
\vbox{\kern7\p@\hbox{.}}\mkern2mu
\raise4\p@\hbox{.}\mkern2mu\raise7\p@\hbox{.}\mkern1mu}}
\makeatother

\title{A note on lower bounds for hypergraph Ramsey numbers}
\author{David Conlon\thanks{St John's College, Cambridge, CB2 1TP, United Kingdom. E-mail: {\tt
D.Conlon@dpmms.cam.ac.uk}}}
\date{}

\begin{document}
\maketitle

\begin{abstract}
We improve upon the lower bound for $3$-colour hypergraph Ramsey
numbers, showing, in the $3$-uniform case, that
\[r_3 (l,l,l) \geq 2^{l^{c \log \log l}}.\]
The old bound, due to Erd\H{o}s and Hajnal, was
\[r_3 (l,l,l) \geq 2^{c l^2 \log^2 l}.\]
\end{abstract}

\section{Introduction}

The hypergraph Ramsey number $r_k (l, l)$ is the smallest number
$n$ such that, in any $2$-colouring of the complete $k$-uniform
hypergraph $K_n^k$, there exists a monochromatic $K_l^k$. That
these numbers exist is exactly the statement of Ramsey's famous
theorem \cite{R30}.

These numbers were studied in detail by Erd\H{o}s and Rado
\cite{ER52}, who showed that
\[r_k (l, l) \leq 2^{2^{\Ddots^{2^{c l}}}},\]
where the tower is of height $k$ and $c$ is a constant that
depends on $k$.

For the lower bound, there is an ingenious construction, due to
Erd\H{o}s and Hajnal (\cite{GRS80}, \cite{EHR65},\cite{EHMR84}),
which allows one to show, for $k \geq 3$, that
\[r_k (l, l) \geq 2^{2^{\Ddots^{2^{c l^2}}}},\]
where this time the tower is of height $k-1$ and $c$ is another
constant depending on $k$. Their construction uses a so-called
stepping-up lemma, which allows one to construct counterexamples
of higher uniformity from ones of lower uniformity, effectively
giving an extra exponential each time we apply it to move up to a
higher uniformity. Unfortunately, it does not allow one to step up
graph counterexamples to $3$-uniform counterexamples, and it is
here that we lose out on the single exponential by which the
towers differ. Instead, we have to start from a different
$3$-uniform counterexample, the simple probabilistic one, which
yields
\[r_3 (l, l) \geq 2^{c l^2},\]
and use that to step up.

Erd\H{o}s was obviously very fond of this problem, offering
$\$500$ for the person who could close the gap between the upper
and the lower bound. As yet, there has been no progress, in the
$2$-colour case, beyond the bounds we have given above. However,
the number of colours seems to matter quite a lot in this problem.
Erd\H{o}s and Hajnal were already aware (again, see \cite{GRS80})
that a variant on their methods could produce a counterexample
showing that indeed
\[r_k (l, l, l, l) \geq 2^{2^{\Ddots^{2^{c l}}}},\]
where now the tower has the correct height $k$ and again $c$
depends on $k$.

Naturally, in the $3$-colour case, one would also expect some
little improvement, and Erd\H{o}s and Hajnal provided just such a
result (unpublished, see \cite{CG98}, though the reader may
consult \cite{EH89} for an earlier attempt), showing that
\[r_3 (l, l, l) \geq 2^{c l^2 \log^2 l}.\]
It is this case that we will look at in this paper, showing that
the bound may be improved rather more substantially to

\begin{theorem}
\[r_3 (l, l, l) \geq 2^{l^{c \log \log l}}.\]
\end{theorem}

Our method is in the stepping-up lemma tradition. It differs,
however, from the lemmas proved in the past in that we make
explicit use of the probabilistic method in our construction. A
rough idea of the proof is that we choose a very dense graph $G$
containing no cliques of size $l$. We then step up to a dense
$3$-uniform graph $H$ and $2$-colour it. The specific form of the
$2$-colouring implies that we cannot contain a monochromatic
$3$-uniform $(l+1)$-clique without $G$ containing an $l$-clique.
The complement of $H$ is coloured with the third colour. It is the
step up of the complement of $G$, and we show, by an involved
argument, that this sparse graph can be chosen in such a way that
the third colour (the step-up of this graph) does not contain an
$(l+1)$-clique. It is this part of the argument which is new and
facilitates our improvement.

Once we have the $3$-uniform case, we can then apply the
stepping-up lemma of Erd\H{o}s and Hajnal, which we state as

\begin{theorem}
If $k \geq 3$ and $r_k (l) \geq n$, then $r_{k+1} (2l + k - 4)
\geq 2^n$.
\end{theorem}

to give the following theorem

\begin{theorem}
\[r_k (l, l, l) \geq 2^{2^{\Ddots^{2^{l^{c \log \log l}}}}},\]
where the tower is of height $k$ and the constant $c$ depends on
$k$.
\end{theorem}

\section{Proof of Theorem 1}

Note that, throughout this section, whenever we use the term
$\log$ we mean log taken to the base 2.

Let $G$ be a graph on $n$ vertices which does not contain a clique
of size $l$ . We are going to consider the complete $3$-uniform
hypergraph on the set
\[T = \{(\gamma_1, \cdots, \gamma_n) : \gamma_i = 0 \mbox{ or }
1\}.\]

If $\epsilon = (\gamma_1, \cdots, \gamma_n)$, $\epsilon' =
(\gamma'_1, \cdots, \gamma'_n)$ and $\epsilon \neq \epsilon'$,
define
\[\delta(\epsilon, \epsilon') = \max\{i : \gamma_i \neq
\gamma'_i\},\] that is, $\delta(\epsilon, \epsilon')$ is the
largest component at which they differ. Given this, we can define
an ordering on $T$, saying that
\[\epsilon < \epsilon' \mbox{ if } \gamma_i = 0, \gamma'_i = 1,\]
\[\epsilon' < \epsilon \mbox{ if } \gamma_i = 1, \gamma'_i = 0.\]
Equivalently, associate to any $\epsilon$ the number $b(\epsilon)
= \sum_{i=1}^n \gamma_i 2^{i-1}$. The ordering then says simply
that $\epsilon < \epsilon'$ iff $b(\epsilon) < b(\epsilon')$.

We will do well to note the following two properties of the
function $\delta$:\\

(a) if $\epsilon_1 < \epsilon_2 < \epsilon_3$, then
$\delta(\epsilon_1, \epsilon_2) \neq \delta(\epsilon_2,
\epsilon_3)$;

(b) if $\epsilon_1 < \epsilon_2 < \cdots < \epsilon_m$, then
$\delta (\epsilon_1, \epsilon_m) = \max_{1 \leq i \leq m-1}
\delta(\epsilon_i, \epsilon_{i+1})$. \\

Now, consider the complete $3$-uniform hypergraph $H$ on the set
$T$. If $\epsilon_1 < \epsilon_2 < \epsilon_3$, let $\delta_1 =
\delta(\epsilon_1, \epsilon_2)$ and $\delta_2 = \delta(\epsilon_2,
\epsilon_3)$. Note that, by property (a) above, $\delta_1$ and
$\delta_2$ are not equal. Colour the edge $\{\epsilon_1,
\epsilon_2, \epsilon_3\}$ as follows:\\

$C_1$, if $\{\delta_1, \delta_2\} \in e(G)$ and $\delta_1 <
\delta_2$;

$C_2$, if $\{\delta_1, \delta_2\} \in e(G)$ and $\delta_1 >
\delta_2$;

$C_3$, if $\{\delta_1, \delta_2\} \not\in e(G)$.\\

Suppose that $C_1$ contains a clique $\{\epsilon_1, \cdots,
\epsilon_{l+1}\}_<$ of size $l+1$. For $1 \leq i \leq l$, let
$\delta_i = \delta(\epsilon_i, \epsilon_{i+1})$. Note that the
$\delta_i$ form a monotonically increasing sequence, that is
$\delta_1 < \delta_2 < \cdots < \delta_l$. Also, note that since,
for any $1 \leq i < j \leq l$, $\{\epsilon_i, \epsilon_{i+1},
\epsilon_{j+1}\} \in C_1$, we have, by property (b) above, that
$\delta(\epsilon_{i+1}, \epsilon_{j+1}) = \delta_j$, and thus
$\{\delta_i, \delta_j\} \in e(G)$. Therefore, the set $\{\delta_1,
\cdots, \delta_l\}$ must form a clique of size $l$ in $G$. But we
have chosen $G$ so as not to contain such a clique, so we have a
contradiction. Similarly, $C_2$ cannot contain a clique of size
$l+1$.

For $C_3$, assume again that we have a monochromatic clique
$\{\epsilon_1, \cdots, \epsilon_{l+1}\}_<$ of size $l+1$, and, for
$1 \leq i \leq l$, let $\delta_i = \delta(\epsilon_i,
\epsilon_{i+1})$. Not only can we no longer guarantee that these
form a monotonic sequence, but we can no longer guarantee that
they are distinct. Suppose, however, that there are $d$ distinct
values of $\Delta$. We will consider the graph $J$ on the vertex
set $\{\Delta_1, \cdots, \Delta_d\}$ with edge set given by all
those $\{\Delta_i, \Delta_j\}$ such that there exists $\epsilon_r
< \epsilon_s < \epsilon_t$ with $\{\Delta_i, \Delta_j\} =
\{\delta(\epsilon_r, \epsilon_s), \delta(\epsilon_s,
\epsilon_t)\}$. Understanding the properties of these graphs is
essential because these graphs are exactly the ones that we will
need to avoid in the complement of $G$ in order to avoid
stepping-up to a complete graph.

How many edges are there in $J$? To begin, note that $\Delta_1$ is
joined to all other $\Delta$, so we have at least $d - 1$ edges.
Suppose that the one occurrence of $\Delta_1$ is at
$\delta_{i_1}$. For the sake of later brevity, note that we may
sometimes refer to these as $\Delta_{1,1}$ and $\delta_{i_{1,1}}$
respectively. Now, let $\Delta_{2,1}$ be the largest $\delta_j$,
at $\delta_{i_{2,1}}$ say, to the left of $\delta_{i_{1,1}}$ (that
is with $j < i_{1,1}$), a region which we will denote by
$R_{2,1}$. Similarly, let $\Delta_{2,2}$ be the largest $\delta$,
occurring at $\delta_{i_{2,2}}$, in the region $R_{2,2}$ which is
to the right of $\delta_{i_{1,1}}$. $\Delta_{2,1}$ (resp.
$\Delta_{2,2}$) must then be joined to every $\delta$ which is to
the left (resp. right) of $\delta_{i_{1,1}}$. Therefore, since
there must be representatives of all remaining $\Delta$ amongst
these $\delta$, we see that between $\Delta_{2,1}$ and
$\Delta_{2,2}$, they must have $d - t_2$ neighbours, where $t_2$
is the number of distinct $\Delta$ amongst $\Delta_{i_{1,1}},
\Delta_{i_{2,1}}, \Delta_{i_{2,2}}$.

Continuing inductively, suppose that we have the collection of
$\Delta_{a, b}$ for all $1 \leq a \leq i-1$ and all $1 \leq b \leq
2^{a-1}$. This collection, consisting of at most $2^{i-1} - 1$ of
the $\delta$, partitions the $\delta$ into at most $2^{i-1}$
regions, which, starting from the left with $\delta_1$ and working
towards $\delta_l$ on the right, we denote by $R_{i,1}, R_{i,2},
\cdots, R_{i, 2^{i-1}}$. Choose, within each region $R_{i,j}$, the
largest $\delta$, which we denote by $\Delta_{i,j}$. Each of these
is necessarily distinct from all $\Delta_{a,b}$ with $1 \leq a
\leq i-1$. Let $t_i$ be the number of distinct $\delta$ given by
the list of numbers $\Delta_{a,b}$ for $1 \leq a \leq i$ and $1
\leq b \leq 2^{a-1}$. Then, since each of the remaining $\Delta$
must lie in one of the regions $R_{i, j}$, we see that at least
one of $\Delta_{i,1}, \cdots, \Delta_{i, 2^{i-1}}$ must be
connected to each of the $d - t_i$ remaining $\Delta$.

We continue this process until we run out of representatives, that
is until the $m$th step, when $t_m = d$. Note that there must be
such an $m$, since we must add at least one new $\Delta$ class at
each step. Note also that $m \geq \log (l + 1) $. This is because,
unless we have used up all of the $\delta$ in our process there
will always be some extra distinct representatives remaining to
consider. So we must have that $2^m - 1$, which is the maximum
number of $\delta$s considered at step $m$, is at least as large
as $l$. Consequently, as $d \geq m$, we also have that $d \geq
\log (l + 1)$.

Now, overall, we have
\[(d - t_1) + \cdots + (d - t_m) = dm - (t_1 + \cdots + t_m)\]
edges. To get a lower bound on this, we need to have upper bounds
for each of the $t_i$. A straightforward upper bound for $t_i$,
following from the fact that $t_i \geq t_{i-1} + 1$, is $t_i \leq
d - m + i$. For small $i$ we can do better, since there we know
that $t_i \leq 2^i - 1$. Therefore, letting $i_0 = \log (d - m + 1
)$, we have
\begin{eqnarray*}
t_1 + \cdots + t_m & = & \sum_{i=1}^{i_0} t_i + \sum_{i = i_0+1}^m
t_i\\ & \leq & 2(d - m + 1) + \sum_{i = i_0+1}^m (d - m + i)\\ & =
& 2(d - m + 1) + \sum_{j = 0}^{m - i_0 - 1} (d - j)\\ & = &
2(d-m+1) + d (m - i_0) - \frac{(m - i_0)(m - i_0 - 1)}{2}.
\end{eqnarray*}
Subtracting this from $d m$, we see that the total number of edges
is at least
\[d i_0 + \frac{(m - i_0)(m - i_0 - 1)}{2} -  2 (d - m + 1).\]
Now, if $d-m+1 \geq \frac{1}{2} \log (l+1)$, we have, since $i_0 =
\log(d-m+1)$, that this is greater than
\[d (\log \log (l+1) - 3).\]
If, on the other hand, $d - m + 1 \leq \frac{1}{2} \log (l+1)$, we
have that $m \geq d + 1 - \frac{1}{2} \log (l+1) \geq d/2 + 1$
(recall that $d \geq \log (l+1)$) and, therefore, the number of
edges is at least
\[\frac{1}{8} (d + 2 - 2 \log \log (l+1))(d - \log \log (l+1))
- 2d \geq \frac{1}{10} d \log (l+1),\] for $l$ large. So, in any
case, for $l$ sufficiently large, we have that the number of edges
is at least
\[\frac{1}{10} d \log \log (l+1).\]

Now, for any graph $J$, let $J'$ be the graph formed by the
process of joining $\Delta_{i, j}$ to all $\Delta$ that have
representatives in the region $R_{i,j}$. If at any stage we find
that we have $\Delta_{i, j_1}$ and $\Delta_{i, j_2}$, both of
which are joined to the same $\Delta$, then we remove one of the
edges arbitrarily, eventually forming a graph $J''$. Every graph
$J$ must contain such a graph. In fact, above, it is the minimum
number of edges in an associated $J''$ that we have counted. The
question we must now ask is, how many distinct $J''$, up to
isomorphism, are there, given that we have a certain $d$?

Consider the set of vertices $V = \{v_1, \cdots, v_d\}$. Choose
the vertex $v_1$ and join it to all other vertices. Now consider
the set $V\char92\{v_1\}$. Up to isomorphism there are at most $d$
different ways to partition this set into two sets $V_{2,1}$ and
$V_{2,2}$, say. Now choose a vertex in each set, say $v_{2,1}$ and
$v_{2,2}$, and join each to all other vertices in their respective
sets. Consider, in turn, the sets $V_{2,1}\char92\{v_{2,1}\}$ and
$V_{2,2}\char92\{v_{2,2}\}$, and partition each into two sets
$V_{3,1}, V_{3,2}$ and $V_{3,3}, V_{3,4}$ respectively. Again, up
to isomorphism there are at most $d$ ways to partition each of the
sets. So, overall, we have at most $d^3$ non-isomorphic classes at
this stage.

Continue in the same way. At the $i-1$st stage, we have sets
$V_{i-1, 1}, \cdots, V_{i-1, 2^{i-2}}$. Choose, in each set
$V_{i-1,j}$, a vertex $v_{i-1,j}$ and join it to every other
vertex in the set. Then partition each set
$V_{i-1,j}\char92\{v_{i-1,j}\}$ into two sets. As always, this can
be done, up to isomorphism in at most $d$ ways. This process stops
when we run out of vertices.

Note that at each step we choose a vertex and then partition an
associated set. Since there are at most $d$ vertices and the
number of ways to partition any set is at most $d$, we conclude
that the number of non-isomorphic graphs $J''$ is at most $d^d$.
(This is, of course, quite a rough estimate, but it is relatively
easy to prove and perfectly sufficient for our purposes.)

We are finally ready to pick the graph $G$. Recall that, for the
first two colours not to contain a $3$-clique of size $l+1$, we
need to choose $G$ so as not to contain a clique of size $l$.
Moreover, for the last colour not to contain a $3$-clique of size
$l+1$, it is sufficient that the complement of $G$, denoted by
$\overline{G}$, does not contain any of the graphs $J''$.

We are going to fix $n = l^{c \log \log l}$, where $c$ is a
constant to be determined, and choose edges with probability $p =
1 - \frac{\log l \log \log l}{l}$. The expected number of cliques
of size $l$ in $G$ is then
\begin{eqnarray*}
p^{\binom{l}{2}} \binom{n}{l} & = & \left(1 - \frac{\log l \log
\log l }{l}\right)^{\binom{l}{2}} l^{c l \log \log l}\\ & \leq &
e^{- \frac{1}{2} l \log l \log \log l} e^{c l \log l \log \log
l}\\ & \leq & e^{- \frac{1}{4} l \log l \log \log l},\\
\end{eqnarray*}
if we take $c \leq 1/4$.

On the other hand, the expected number of graphs $J''$ of order
$d$ that we can expect to find in $\overline{G}$ is at most
\begin{eqnarray*}
d^d (1-p)^{\frac{1}{10} d \log \log (l+1)} n^d & \leq &
\left(\frac{\log l \log \log l}{l}\right)^{\frac{1}{10} d \log
\log (l+1)} (d n)^d\\ & \leq & e^{- \frac{1}{20} d \log l \log
\log (l+1)} l^{2 c d \log \log l}\\ & \leq & e^{- \frac{1}{40} d
\log l \log \log l},
\end{eqnarray*}
if we take $c \leq 1/80$ and $l$ sufficiently large.

Adding over the expected number of cliques in $G$ and the expected
number of copies of graphs $J''$ in $\overline{G}$ for all $l$
possible values of $d$, we find that, for $l$ sufficiently large,
the expected value of all such graphs is less than one. We can
therefore choose our graph $G$ in such a way that it does not
itself contain a clique of size $l$ and its complement
$\overline{G}$ does not contain any of the graphs $J''$. The
result follows.

\end{document}